\documentstyle[]{article}
\author {Anna Erschler (Dyubina)}
\date{  e-mail: annadi@post.tau.ac.il, 
erschler@pdmi.ras.ru}
\title{On the asymptotics of drift}

\begin{document}
\maketitle

\newtheorem{rem}{Remark}
\newtheorem{definition}{Definition}
\newtheorem{lemma}{Lemma}
\newtheorem{corollary}{Corollary}
\newtheorem{proposition}{Proposition}
\newtheorem{theorem}{Theorem}
\newcommand{\erw}{{\mathop{\rm E}}}
\newcommand{\Z}{\rm{\hbox to 9pt{Z\hss Z}}}
\newcommand{\R}{{\rm R}}
\newcommand{\tgot}{{\sf t}}
\newcommand{\zz}{\Z/2\Z}

\section{Introduction}
In this paper we consider such symmetric random walks on finitely
generated groups  that the support of the corresponding measure
is finite and generates the group.
To any such walk given by a probability measure
 $\mu$ we can assign a function
 (see \cite{Ve}, \cite{Ve2})
$$
L_n=\erw_{\mu^{*n}} l(g),
$$
where $l$ denotes the word length corresponding to a fixed system
of generators
  and
$\mu^{*n}$ is $n$-th convolution of $\mu$.
The function $L_n$ is called the drift (or escape) of the random walk.
The function $L_n$ measures how fast in average the random walk
is moving away from the
origin.
This function carries much information.
For example, it is asymptotically linear iff
the entropy of the random walk is positive
(see \cite{Ve2}, \cite{Va}). In particular this function has linear
asymptotics for any nonamenable group.
In many other examples (for example, for any
Abelian group)
 $L_n$ is asymptotically equivalent to $\sqrt{n}$. The first example of
a non-trivial
(that is nonlinear and not $\sqrt{n}$)
asymptotics of
 $L_n$
is constructed in \cite{Ja}. There I prove that
for some simple random walk on a group
$L_n$ is asymptotically $n/\ln(n)$. In \cite{Ja}
some other examples were anounced.
This paper contains a proof for these assertions.
We say that a random walk is simple if the corresponding measure
is equidistributed on some finite symmetric generating set.
In this paper we prove that for any positive integer $k$
there exist a group and a simple random walk on it such that
$$
L_n \asymp n^{1-\frac{1}{2^k}}.
$$

\section{Auxiliary statements about random walk on a line}

For the proof of the main result we need some standard lemmas about
symmetric random walks on a line.

{\bf Definition.}
The  {\it range} of the random walk is the number of
different elements of the group visited during the first
$n$ steps ($n\ge 0$).
Let $R_n$ denote the range.
\begin{lemma} Consider a simple symmetric random walk on
$\Z$. Then
$$ \erw[R_n]\asymp \sqrt{n}$$
\end{lemma}
{\bf Proof.}
See for example \cite{Sp}

As before we consider simple symmetric random walk on
 $\Z$. Let $b_i^{(n)}$ be the number of times that the random
walk visited the element $i$ during the first
 $n$ steps.
 (In this case $b_i^{(n)}$ is called the
{\it local time} of the random walk).
\begin{lemma} There exist $C >0$ and
 $p>0$ such that for any $n$

$$
\Pr \left[
\max_{i} b_i^{(n)}\le C\sqrt{n} \right] >p
$$
\end{lemma}
{\bf Proof.}
Normalized local time of the random walk tends uniformly to that
of the Brownian motion.
In fact,
 \cite[ Theorem 3.2,  Chapter 5]{IB} states that
$$
\lim_{n \to \infty}{\Pr \left(\sup_{(t,x)\in [0,T]\times \R}
|\tgot_n(t,x)-\tgot(t,x)|>\varepsilon \right)}=0,
$$
In particular
$$
\lim_{n \to \infty}
{\Pr \left(\sup_{x \in \R} |\tgot_n(1,x)-\tgot(1,x)|
> \varepsilon \right)}=0.
$$
Here $\tgot(t,x)$ is the local time of the Brownian motion
and $\sqrt{n} \tgot_n(t,x) $ is the number of times
 the random walk has visited the point
$[x\sqrt{n}=i]$ up to the moment $[nt]$.
Therefore in the preceding notations $\tgot_n(t,x)$ is equal to
$b_i^{(n)}/\sqrt{n}$.
For the distribution of the supremum of the local time see
 \cite{Bor} .
It suffices to know that there exists
 $C'>0$ such that
$\sup_{x} \tgot(1,x) \le C'$
with a positive probability not depending on $n$.

\begin{proposition}
Let $0<\alpha<1$. For simple symmetric random walk on $\Z$
$$\erw \left[\sum_{i} \left(b_i^{(n)} \right)
^\alpha\right] \asymp n^{\frac{\alpha+1}{2}}.$$
\end{proposition}
{\bf Proof.}
Note that
$R_n=\# \{i: b_i^{(n)}\ne 0\}$ É $\sum b_i^{(n)} =n$.

Note also that since $x^\alpha$ is concave we have
$$
 \sum_{i} \left(b_i^{(n)} \right)
^\alpha \le R_n \left(\frac{n}{R_n}\right)^\alpha
=R_n^{1-\alpha}n^\alpha.
$$
Then, using the fact that $x^{1-\alpha}$ is also concave we get
$$
\erw \left[\sum_{i} \left(b_i^{(n)} \right)
^\alpha \right]
\le \erw \left[R_n^{1-\alpha}n^\alpha\right] =
 n^\alpha\erw \left[R_n^{1-\alpha}\right]\le
 n^\alpha\erw[R_n]^{1-\alpha} \asymp
n^{\alpha+(1-\alpha)/2}.
$$

The last passage follows from that, as it proved in lemma 1,

$\erw[R_n] \asymp \sqrt{n}$.

Now let us estimate
$\erw \left[\sum_{i} (b_i^{(n)})^\alpha \right]$
 from below.
By lemma 1
$$
\Pr \left[ \max_{i} b_i^{(n)}\le C\sqrt{n} \right] > p.
$$
Consider
$\min {\sum_{i} (b_i^{(n)})^\alpha}$
 assuming that $b_i^{(n)}\ge 0$,
$$
\max_{i} b_i^{(n)}\le C\sqrt{n}
$$
and
$$
\sum_{i} b_i^{(n)}=n.
$$

If $C\sqrt{n}\le n$ then the maximum is achieved when all
 $b_i^{(n)}$ but one are equal to
 $C\sqrt{n}$,
since $x^\alpha$ is concave.

 Therefore
$$
\sum_i (b_i^{(n)})^\alpha \ge [n/(C\sqrt{n})] (C\sqrt{n})^\alpha.
$$
Hence there exists $K$ depending only on $C$ such that
$$
\sum_{i} (b_i^{(n)})^\alpha\ge K n^\frac{1+\alpha}{2}.
$$
Therefore
$$
\erw \left[\sum_{i} (b_i^{(n)})^\alpha \right]\ge p K n^\frac{1+\alpha}{2}.
$$
This completes the proof of the proposition.

\section{Examples of drift}

To construct the examples we need the definition of wreath product.

\begin{definition}
Wreath product of $A$ and $B$ is a semidirect product
$A$ and $\sum_A B$,  where $A$ acts on
$\sum_A B$ by shifts:
if $a\in A$, $f:A \to B, f\in \sum_A B$,  then
$f^a(x)=f(xa^{-1}), x\in A$.
Denote the wreath product by $A\wr B$.
\end{definition}
Recall that the wreath product of finitely generated groups is also
finitely generated.
Random walks on wreath products were used for construction
of many non-trivial examples
 (see \cite{VK}, and
also \cite{Ja}, \cite{Ja2}).

The following lemma is the main tool in the proof of the main
result.

\begin{lemma}
Let $a_1, a_2,..., a_k$ be generators of $A$.
Consider the measure uniformely distributed on these
generators and their inverses.

Suppose that
$$L_n^A\asymp n^\alpha.$$
Then for some simple random walk on  $B=\Z\wr A$ we have
$$L_n^B\asymp n^{\frac{\alpha+1}{2}}.$$
\end{lemma}
{\bf Proof.}
Let us make the following notation. Let $a\in A$.
Let $e$ be the identity element. Then
$\tilde{a}^e$ denotes the function from $\Z$ to $A$ such that
$\tilde{a}^e(0)=a$ and  $\tilde{a}^e(x)=e$ for any $x\ne 0$. Let $a^e$
be equal to
$(e,\tilde{a}^e)$,  and $\delta=(1,e)$.

Consider the following set of generators of $B$:
$$(a_j^e)^p \delta (a_k^e)^q,$$
where $p=0,1$ or $-1$,
 $q=0,1$ or $-1$ and
$1\le j,k \le n$.

Consider the simple symmetric random walk on
 $B$ corresponding to this set of generators.

Suppose that at the moment $n$ this random walk hits $(i,f)$
($i \in \Z$, $f: \Z \to Z$). For each $j$
put
$c_j=l_A(f(j))$.
Then
 $$ \frac{1}{2} \sum_{j} c_j \le l(i,f)\le 2(\sum_{j} c_j+R)$$

Multiply $(i,f)$ by one of the generators. It changes the value of
$f$ at most at $2$ points. At each point the value can be shifted by a
generator of $A$. This proves the first inequality.

At most two multiplications by generators of $B$ suffice to move
$(i,f)$
into $(i, f')$, $f'$ conincides with $f$ everywhere outside $i$,
$f(i) \ne f'(i)$. From this follows the second inequality.

Hence
$$
\frac{1}{2}\erw \left[\sum c_j \right]
=\frac{1}{2}\sum\erw[c_j] \le \erw [l(i,f)] \le 2(\sum\erw[c_j]+\erw[R])
$$

The projection of this random walk on $\Z$ is symmetric and simple.
Note that in each $j$ there is some random walk on $A$.
The measure that defines this random walk is concentrated on
 $a_1, a_2,... a_k$, their inverses, and
 $e$.
Note that it is uniformely distributed on
  $a_1, a_2,...a_k$ and their inverses.
Let $\tilde{L}_n^A$ be the drift corresponding to this new measure.
It is clear that
$$\tilde{L}_n^A \asymp L_n^A \asymp n^\alpha.$$

 We can think that we get
$f(j)$ as the position of the  random walk on $A$ after $2b_j^{(n)},
 2b_j^{(n)}-1$
or $2b_j^{(n)}-2$ steps.
The number of steps depends on the fact whether
 $i$ is equal to $0$ or to the final
 (after $n$ steps on $B$) point $i$of the random walk on $\Z$.
 Therefore

$$\erw[c_i | b_i^{(n)}=b, j\ne 0, j\ne i]=\tilde{L}_{2b}^A$$
$$\erw[c_i | b_i^{(n)}=b, j=0, j\ne i]=\tilde{L}_{2b-1}^A$$
$$\erw[c_i | b_i^{(n)}=b, j\ne 0, j=i]=\tilde{L}_{2b-1}^A$$
$$\erw[c_i | b_i^{(n)}=b, j=0, j=i]=\tilde{L}_{2b-2}^A$$

Since $0<\tilde{L}_n^A\asymp n^\alpha$, there exist $C_2,C_3>0$ such that
for any
 $b_j^{(n)}$ and $m=2b_j^{(n)}-2, 2b_j^{(n)}-1, 2b_j^{(n)}$

$$
C_2 (b_j^{(n)})^\alpha\le \tilde{L}_{m}^A \le C_3 (b_j^{(n)})^\alpha
$$

Then
$$
C_2 \erw[(b_j^{(n)})^\alpha] \le \erw[c_j] \le C_3 \erw[(b_j^{(n)})^\alpha]
$$
and hence
$$
\frac{1}{2}C_2 \erw \left[\sum (b_i^{(n)})^\alpha \right] \le
\erw[l(a,f)]=L_n^B \le 2(C_3
         \erw \left[\sum (b_i^{(n)})^\alpha \right]+\erw[R])
$$
By proposition 1
$$
\erw\left[\sum (b_j^{(n)})^\alpha\right]\asymp n^{(1+\alpha)/2},
$$
and
$$
\erw[R]\asymp \sqrt{n}\le   n^{(1+\alpha)/2}.
$$
Hence
$$
L_n^B\asymp    n^{(1+\alpha)/2}.
$$
This completes the proof of the lemma.

As a consequence of lemma 3 we get the following theorem.

\begin{theorem}
Consider the following groups that are defined recurrently
$$G_1=\Z; G_{i+1}=\Z\wr G_i$$
Then for some simple random walk on
 $G_i$ we have
$$
L_n^{G_i} \asymp n^{1-\frac{1}{2^i}}
$$
\end{theorem}

{\bf Proof.} We prove the theorem by induction on
 $i$.
The base is $i=1$. In this case $G_i=\Z$ and
$L_n$ is asymptotically equivalent to $\sqrt{n}$
(see for example \cite{Sp}).
Induction step follows from the previous lemma.

The author expresses her gratitude to A.M.Vershik for stating
 the problem and many helpful discussions.
I am also grateful to M.I.Gordin,
I.A.Ibragimov, and
A.N.Borodin for useful advices.

\end{document}